\documentclass[
12pt]{article}

\usepackage{geometry}
\usepackage{color}
\usepackage{amssymb}

\title{{\normalsize\tt\hfill\jobname.tex}\\
Yet again on polynomial convergence 
for SDEs with a gradient-type drift
\author{A.Yu.~Uglov\footnote{National Research University Higher School of Economics, Russia; email: auglov @ hse.ru},     A.Yu.~Veretennikov\footnote{University of Leeds, United Kingdom; National Research University Higher School of Economics, \& Institute for Information Transmission Problems, Russia; email: a.veretennikov @ leeds.ac.uk} 
}
}

\date{}
\begin{document}
\maketitle
\newtheorem{Theorem}{Theorem}
\newtheorem{Lemma}{Lemma}
\newcommand{\eps}{\varepsilon}
\newtheorem{Proposition}{Proposition}
\newtheorem{remark}{Remark}

\begin{abstract}
Bounds on convergence rate to the invariant distribution for a class of stochastic differential equations (SDEs) 
are studied.

\end{abstract}

\section{Introduction}

Let us consider a stochastic differential equation in $R^d$

\begin{equation}\label{eq1}                                  
dX_t = dB_t -  \nabla U (X_t)\,dt
\end{equation}
with initial data
\begin{equation}\label{eq2}                                   
X_0 = x.
\end{equation}
Here $B_t, \,\, t\geq 0$ is a $d$-dimensional Brownian motion, $X_t$
takes values in $R^d$, $U$ is a
non-negative function, $U(0)=0$ and $\lim_{|x|\to\infty}U(x) = +\infty$.
Function $U$ is assumed to be locally bounded and locally $C^1$. The aim of this paper is to establish ergodic properties of the Markov process $X_t$, namely,
existence and uniqueness of its invariant probability measure, and to estimate 
convergence rate to the invariant measure which rate bound would not depend on the first derivatives of the function $U$.
Such a problem -- about bounds not depending explicitly on \(\nabla U\) --  was posed and in some particular case solved in \cite{Ver01}. Here we extend and relax some of the assumptions from \cite{Ver01}.  
It is widely known that the rate of convergence may be derived from  the estimates of the type

\begin{equation}\label{eq3}                                
\mathbb E_x \tau^k \leq C(1+|x|^m),
\end{equation}

\begin{equation}\label{eq4}                              
\sup_{t\ge 0} \mathbb E_x |X_t|^m \leq C(1+|x|^{m'}),
\end{equation}
for some $k>1$, $m, m', C > 0$, where $\tau = \inf (t\geq 0: \,\,
|X_t| \le K)$ for some $K>0$, see, e.g.,  \cite{Ka, Ve1}, et al. In
particular, for SDEs (\ref{eq1}) with a  {\em bounded} $\nabla U$
it can be derived from (\ref{eq3}) and (\ref{eq4}) that
\begin{equation}\label{eq5}                                  var(\mu^x_t-\mu^{inv}) \leq  P(x)(1+t)^{-k'},
\end{equation}
with any $k'<k$ and some function $P$ growing in \(x\) at infinity.

~

\noindent
The bounds like (\ref{eq3}) under various assimptions were obtained for various classes of
processes by many authors, see, in particular,  \cite{AIM, Ka, La},
\cite{MW} -- \cite{Ve2} and the references therein; yet, for SDEs all assumptions were usually -- except the paper \cite{Ver01} -- stated in terms of \(\nabla U\). See also \cite{GRS, Mal}
where stronger sub-exponential bounds were established under another standing assumption. In
\cite{Ve1} and \cite{Ve2} a recurrence condition
$$
-p = \limsup (\nabla U(x),x) < 0
$$
was used to get bounds like  (\ref{eq5}). Here the goal is to use some analogue of the latter condition but in terms of the limiting behavour of the function \(U\) itself, similar to \cite{Ver01} but under weaker assumptions.

\section{Main results}
\subsection{Earlier results}
Recall briefly some earlier results from \cite{Ver01} where, in fact, a little more general equation was considered.
Assume 
\begin{equation}\label{eq7a}                          
\sup_{x,x':\, |x-x'|\le 1}\; (U(x)-U(x')) < \infty
\end{equation}
and let the structure of the function $U$ be  as follows:

\begin{eqnarray}\label{eq7}                                 
U(x) = U^1(x) + U^2(x),
\quad U^1(x) = V(|x|), \quad <U^2(x),x> \equiv 0.
\end{eqnarray}
The function \(V\) here is assumed in the class \(C^1(0, \infty)\). 
In particular, the "essential" divergent part $U^1$ of the drift has a
central symmetry property while another divergent part $U^2$ is
orthogonal to the direction $x$ at any point $x$.
Let the following recurrence condition is satisfied,

\begin{equation}\label{eq8}                                  
-\lim_{ \xi \to\infty} \frac{V(\xi)}{
\log \xi} + d = -p < 0.
\end{equation}

\begin{Proposition}[\cite{Ver01}]\label{Theorem1}
Let (\ref{eq7a})--(\ref{eq7}) and (\ref{eq8})  with $p>1/2$ be satisfied. Then for any $0<k<p+1/2$ and 
$\varepsilon>0$ small enough the estimate (\ref{eq3}) holds with 
$m=2k+\varepsilon$ and some  $C=C_\varepsilon$ and the estimate 
(\ref{eq4}) is valid with any $m<2p-1$ and $m'=m+2\varepsilon$. 
Moreover, there exists a unique invariant measure for the Markov process 
$X_t$. 
\end{Proposition}

\begin{Proposition}[\cite{Ver01}]\label{Theorem2}
Let (\ref{eq7a})--(\ref{eq7}) and (\ref{eq8})  with $p>1/2$ be satisfied. Then the bound (\ref{eq5}) holds true with any $k'<k<p+1/2$ and
$\tilde P(x) = C_\varepsilon (1+|x|^{m})$, $m=2k+\varepsilon$
with any $\varepsilon>0$ small enough and some $C$.  
If, moreover, $p>3/2$ then the 
bound (\ref{eq5a}) holds true with any $k'<k<p-1/2$ and
$\hat P(x) = C_\varepsilon (1+|x|^{m})$, $m=2k+\varepsilon$
with any $\varepsilon>0$ small enough and some $C$.
\end{Proposition}
The assumption $p>3/2$ relates to the critical value $3/2$ in \cite{Ve1}.

\subsection{\bf New results}
Below \([a]\) denotes the integer value of \(a\in \mathbb R^1\). 
\begin{Theorem}\label{Thm1}
Let there exist \(1/2<p_2\le p_1\) such that  
\begin{equation}\label{eq8b}
0<p_2 \le \frac{V(\xi)}{\log \xi} - d \le p_1, 
\end{equation}
 for all \(\xi >0\) which are large enough by the absolute value. Then, 
the bound (\ref{eq4}) holds true with \(m' = m + 2(p_1 - p_2)\)  and \(m=2k(1+p_1-p_2)\). Moreover, 
for any  positive integer value of \(\displaystyle k<1+\frac{2p_2-1}{2(1+p_1-p_2)}\) and \(m=2k(1+p_1-p_2)\), the bound (\ref{eq3}) holds. Moreover 
there is a unique invariant probability measure \(\mu^{inv}\), and for any \(0<k'<k\), and for any \(t\ge 0\), 
\begin{equation}\label{eq5a}
var(\mu^x_t-\mu^{inv}) \leq P(x)(1+t)^{-k'}, 
\end{equation}
with some polynomial function \(P(x)\).
\end{Theorem}

\begin{remark}
Note that \(k=1\) is included in the range of values for which the bound (\ref{eq3}) will be established. 
The assumption (\ref{eq8b}) may be replaced by a similar one with \(\limsup_{|\xi|\to\infty}\) and \(\liminf_{|\xi|\to\infty}\) instead of exact inequalities which may or may not change slightly the resulting statement depending on whether or not  the value \(\displaystyle \frac{2p_2-1}{2(1+p_1-p_2)}\) is integer. Also, depending on whether the same value is integer, the range of \(k\) for which the bound (\ref{eq3}) holds true may change a bit. We do not pursue the inspection of all these possible  changes here. Let us mention that the assumption (\ref{eq7a}) is needed for the ``local mixing'' which explaination may be read in \cite{Ver01} in  detail. 
\end{remark}

\section{Proof}

{\bf 1.} 
As in \cite{Ver01}, due to comparison theorems for SDEs with reflection and
the assumption on the structure of the drift one gets,

\begin{eqnarray}
& \displaystyle  |X_t| \le y_t,
 \nonumber \\ \nonumber  \\
& \displaystyle dy_t = d\bar w_t + \left(\frac{d}{y_t} - V'(y_t) \right)dt + d\varphi_t
\equiv d\bar w_t - \bar V'(y_t)dt + d\varphi_t,\label{eqy}
\end{eqnarray}
where $\bar w$ is a 1-dimensional Wiener process, $y$ is a
solution of the SDEs above with a non-sticky boundary condition
at (any) point $K>0$, $\varphi$ is its local time at $K$, $\bar
V'(y) = V'(y) - d/y $; in other words, we let 
\[
\bar V(y) = V(y) - d \ln y, \quad y>0.
\] 
Condition (\ref{eq8b}) can be
rewritten in the form
\[
\xi^{2p_2} \le \exp(2\bar V(\xi)) \le \xi^{2p_1}, \quad \xi \ge K.
\]

~

\noindent
{\bf 2.}
The invariant density of the process $\xi_t$ with $K=|x|$ has a form
$$
C(|x|) \exp\left(-2\bar V(y)\right), \quad y > |x|.
$$
The normalizing identity implies the estimation from above (under \(2p_2 > 1\)),
$$
C(|x|) = \left(\int_{|x|}^\infty \, \exp(-2\bar V(y)\,dy)\right)^{-1}
\le \left(\int_{|x|}^\infty \xi^{-2p_2}\,dy)\right)^{-1}
= (2p_2-1) |x|^{2p_2  -1}, 
$$
for the values of \(|x|\) large enough. For smaller values of \(|x|\), convergence of the integral cannot be destroyed because in some bounded neighbourhood of zero the function \( \exp\left(-2\bar V(y)\right)\) is bounded. Note that for small values of \(|x|\) the expressions \(\displaystyle \left(\int_{|x|}^\infty \, \exp(-2\bar V(y)\,dy)\right)^{-1}\) are smaller, which means that in all cases for some \(C_0\), 
\[
C(|x|) \le (2p_2-1) |x|^{2p_2  -1} \wedge C_0. 
\]

~

\noindent
{\bf 3.} The inequality  (\ref{eq4}) with any real value  
\(m< 2p_2 -1\) and with \(m' = m + 2(p_1 - p_2)\) (where \(m'\) may not be necessarily integer either) 
follows from a direct calculation,
\begin{eqnarray*}
& \displaystyle \mathbb E_x |X_t|^m \le \mathbb E_{|x|} |y_t|^{m} 
 \le C(|x|)\,\int_{|x|}^{\infty} \,\xi^m 
\exp(-2\bar V(\xi))\,d\xi 
 \\\\
& \displaystyle \le (C |x|^{2p_1-1} \wedge C_0)\,\int_{|x|}^{\infty} \,\xi^m \xi^{-2p_2}\,d\xi \le C |x|^{ m + 2(p_1-p_2)}
\end{eqnarray*}
(here the constants \(C\) may be different on different lines and even on the same  line),
which is true for any $x$ large enough, due to comparison theorems for the processes $y_t$
with different initial data $y_0$. For {\em any \(x\)} -- not necessarily small - this implies the bound (\ref{eq4}), as required. 

~

\noindent
{\bf 4.}
Denote $v^q(\xi) = \mathbb E_\xi\gamma^q$ for any integer \(q\ge 0\), $\gamma = \inf (t:\; y_t\le K)$ and let
$L$ denote the generator of $y_t$. By virtue of the identity
$$
\left( \int_0^\gamma 1\,dt \right)^q = q \int_0^\gamma\,
\left(\int \, 1\,ds\right)^{q-1}\,dt,
$$
it follows,
$$
v^q(y_0) = q \mathbb E_{y_0} \int_0^\gamma v^{q-1}(y_t)\,dt,
$$
for any $q$ such that the integral in the right hand side converges. In turn, this implies an equation (for example, by It\^o's  or Dynkin's  formula) 
\begin{equation}\label{eq10}                                  
Lv^{q} = - q v^{q-1}, \quad (q\ge 1)
\end{equation}
(cf. with \cite{D} theorem 13.17 where the equation is explained differently and under another stronger assumption). Evidently, one
boundary value for the latter equation is $v^q(K)= 0$.  Concerning  the ``second boundary value'' usual for a PDE of the second order, it is seemingly missing
here. The justification of the formula for solution below can be done by the following limiting procedure. Let \(N>K\) be the second boundary (later on \(N\) would go to infinity). Let 
$v_N^q(\xi) = \mathbb E_\xi\gamma_N^q$ for any integer \(q\ge 0\), $\gamma_N = \inf (t:\; y^N_t\le K)$, where the process \(y^N_t\) is a solution of the equation similar to (\ref{eqy}) but with another non-sticky reflection at \(N\). Note that all solutions are strong and, hence, may be constructed on the same probability space; see, e.g., \cite{Ver81} for SDEs with one boundary, and results from this paper are easily extended on the case with two finite boundaries. Apparently, \(y^N_t\le y_t\) for any \(t\) and \(N\), and \(\gamma_N\uparrow \gamma\) as \(N\uparrow \infty\).  So, by the monotone convergence, \(v_N^q \uparrow v^q\) for all values of \(q\) (even if the limit \(v^q\) is not finite). Then the sequence of the functions \(v_N^q(\xi)\) satisfies the equations (\ref{eq10}) with boundary conditions 
\[
v_N^q(K) =0, \quad (v_N^q)'(N) = 0. 
\] 
The formula for solution of such an equation reads, 
\[
v_N^q(\xi) = 2q \int_K^\xi\, \exp(2\bar V(y_1))\,dy_1\,
\int_{y_{1}}^N v_N^{q-1}(y_2) \exp(-2\bar V(y_2))\,dy_2, \quad K\le \xi \le N,
\]
which may be verified by a direct calculation. 
Hence, by induction, the function \(v^q(\xi)\) is given by the formula  via the function $v^{q-1}$,
\begin{equation}\label{eq11}                                  
v^q(\xi) = 2q \int_K^\xi\, \exp(2\bar V(y_1))\,dy_1\,
\int_{y_{1}}^\infty v^{q-1}(y_2) \exp(-2\bar V(y_2))\,dy_2.
\end{equation}
By another induction this implies the inequalities (assuming \(v^0 \equiv 1\)):
\begin{eqnarray*}
v^1(\xi) =2\int_K^\xi\, \exp(2\bar V(y_1))\,dy_1\,
\int_{y_{1}}^\infty v^{0}(y_2) \exp(-2\bar V(y_2))\,dy_2
 \\\\
=  2\int_K^\xi\, \exp(2\bar V(y_1))\,dy_1\,
\int_{y_{1}}^\infty  \exp(-2\bar V(y_2))\,dy_2
\le  2\int_K^\xi\, y_1^{2p_1}\,dy_1\,
\int_{y_{1}}^\infty  y_2^{-2p_2}\,dy_2
 \\\\
= C\int_K^\xi\, y_1^{2p_1-2p_2+1}\,dy_1
= C (\xi^{2(p_1-p_2)+2} - K^{2(p_1-p_2)+2}) 
\le C \xi^{2(p_1-p_2)+2},  
\end{eqnarray*}
under the condition that \(p_2>1/2\) (otherwise the inner integral diverges).
Further, 
\begin{eqnarray*}
& \displaystyle v^2(\xi) =4\int_K^\xi\, \exp(2\bar V(y_1))\,dy_1\,
\int_{y_{1}}^\infty v^{1}(y_2) \exp(-2\bar V(y_2))\,dy_2
 \\\\
& \displaystyle \le C \int_K^\xi\, \exp(2\bar V(y_1))\,dy_1\,
\int_{y_{1}}^\infty  y_2^{2(p_1-p_2)+2} \exp(-2\bar V(y_2))\,dy_2
 \\\\
& \displaystyle \le C  \int_K^\xi\, y_1^{2p_1}\,dy_1\,
\int_{y_{1}}^\infty  y_2^{2(p_1-p_2)+2 - 2p_2}\,dy_2
 \\\\
& \displaystyle = C  \int_K^\xi\, y_1^{2p_1}\,dy_1\,
y_{1}^{2p_1-4p_2+3}
= C (\xi^{4(p_1-p_2)+4} - K^{4(p_1-p_2)+4})
 \\\\
& \displaystyle \le C \xi^{4(p_1-p_2+1)},  
\end{eqnarray*}
where in the calculus it was assumed that 
\(2p_1 - 4p_2 + 2 < -1\), that is, that \(p_1 < 2p_2 -3/2\), otherwise the inner integral in the calculus diverges. Since from the beginning \(p_1\ge p_2\), for the value of \(p_2\) this means that compulsory \(p_2 > 3/2\). 

~

\noindent
Next, 
\begin{eqnarray*}
& \displaystyle v^3(\xi) =6\int_K^\xi\, \exp(2\bar V(y_1))\,dy_1\,
\int_{y_{1}}^\infty v^{2}(y_2) \exp(-2\bar V(y_2))\,dy_2
 \\\\
& \displaystyle \le C \int_K^\xi\, \exp(2\bar V(y_1))\,dy_1\,
\int_{y_{1}}^\infty  y_2^{4(p_1-p_2+1)} \exp(-2\bar V(y_2))\,dy_2
 \\\\
& \displaystyle \le C  \int_K^\xi\, y_1^{2p_1}\,dy_1\,
\int_{y_{1}}^\infty  y_2^{4(p_1-p_2+1) - 2p_2}\,dy_2
 \\\\
& \displaystyle = C  \int_K^\xi\, y_1^{2p_1}\,dy_1\,
y_{1}^{4p_1-6p_2+5}
= C (\xi^{6(p_1-p_2+1)} - K^{6(p_1-p_2+1)})
\le C \xi^{6(p_1-p_2+1)}. 
\end{eqnarray*}
For the inner integral to converge, the values of \(p_1, p_2\) must satisfy \(4p_1 - 6p_2 + 4 < -1\), that is,  \(\displaystyle p_1 < \frac32 p_2-\frac54\). Due to the condition \(p_1 \ge p_2\), for \(p_2\) this compulsory implies \(\displaystyle p_2 > \frac52\). 
Note that, as usual, constants \(C\) may be different for any \(q\) and even from line to line. It looks plausible that the general formula -- as long as the integrals converge -- reads, 
\begin{equation}\label{vq}
v^{q} (\xi) \leq C_{q} \xi^{2q(1+p_1-p_2)}. 
\end{equation}
The base being already established, let us show the induction step. 
Assume that for $q=n-1$ the formula is valid with some constant \(C_{n-1}\), that is,  
\[
v^{n-1} (\xi) \leq C_{n-1} \xi^{2(n-1)(1+p_1-p_2)}.
\]
Then for $q=n$ (as long as the integrals in the calculus below converge) we have, 
\begin{eqnarray*}
& \displaystyle v^{n}(\xi) =2n \int_K^\xi\, \exp(2\bar V(y_1))\,dy_1\,
\int_{y_{1}}^\infty v^{n-1}(y_2) \exp(-2\bar V(y_2))\,dy_2
 \\\\
& \displaystyle \le 2n \int_K^\xi\, {y_1}^{2p_1}\,dy_1\,
\int_{y_{1}}^\infty C_{n-1}\, y_2^{2(n-1)(p_1-p_2+1)} {y_2}^{-2p_2}\,dy_2
 \\\\
& \displaystyle = {C_n}^{}  \int_K^\xi\, y_1^{2p_1}\,y_1^{2n-1+2(n-1)p_1-2np_2}dy_1= {C_n}^{}  \int_K^\xi\, y_1^{2n-1+2np_1-2np_2}dy_1
\le C_n \xi^{2n(p_1-p_2+1)}. 
\end{eqnarray*}
Hence, indeed, by induction the formula (\ref{vq}) is established. The values of \(q\) for which the integrals in the calculus converge must satisfy the bound
\[
2(q-1)(1+p_1-p_2) -2p_2 < -1, 
\]
that is, 
\[
q < q_0:= 1 + \frac{2p_2-1}{2(1+p_1-p_2)} = \frac{1+2p_1}{2(1+p_1-p_2)}. 
\]
As a consequence, it is compulsory that \(p_2 > q-1/2\). Recall that in this paper only integer values of \(q\) are used; however, \(q_0\) introduced above may not be necessarily integer, but in any case \(q_0 >  1\). 
Also, note that if \(p_1=p_2 = p\) as in \cite{Ver01}, then the latter inequality \(q<q_0\)  reduces to \(q< p + 1/2\),  precisely as in \cite{Ver01}.

~

{\bf 5.}
By virtue of the established bounds (\ref{eq3})--(\ref{eq4}), the bound (\ref{eq5}) on convergence towards the stationary measure follows from various sources (cf., e.g.,  \cite{Ve2, Ver01}, et al.) and, hence, in this brief presentation we skip the details of this step. 
The existence of the invariant probability 
measure may be justified via 
the Harris--Khasminsky principle based on 
(\ref{eq3}) with any \(k\ge 1\). Its uniqueness follows, for example, from the bound (\ref{eq4}). This completes the proof of the Theorem \ref{Thm1}. 

\section*{Acknowledgements}
For the first author 
this study has been funded by the 
Russian Foundation for Basic Research grant \mbox{17-01-00633$\_$a}.
For the second author
this study has been funded by the Russian Academic Excellence Project '5-100'
and  
by the Russian Science Foundation project 17-11-01098.


\begin{thebibliography}{99}


\bibitem{AIM} Aspandiiarov, S., Iasnogorodski, R., Menshikov, M. Passage-time
moments for nonnegative stochastic processes and an application to reflected
random walks in a quadrant, Ann. Probab. 1996, 24, 2, 932-960.

\bibitem{D} Dynkin, E.B. Markov processes. New York, Academic Press, 1965.

\bibitem{GRS} Ganidis, H., Roynette, B., Simonot, F. Convergence rate of some
semi-groups to their invariant probability, Stoch. Proc. Appl., 1999, 79(2),
243-263.


\bibitem{Ka} Kalashnikov, V.V. The property of $\gamma$-reflexivity for Markov
sequences. (English. Russian original) Sov. Math., Dokl. 1973, 14, 1869-1873;
translation from Dokl. Akad. Nauk SSSR, 1973, 213, 1243-1246.


\bibitem{La} Lamperti, J. Criteria for stochastic processes II: passage-time
moments. J. Math. Anal. and Appl., 1963, 7, 127-145.

\bibitem{Mal} Malyshkin, M.N. On sub-exponential mixing and convergence rate
for diffusion processes, Toeirya Veroyatn. i ee Primenen., 2000.

\bibitem{MW} Menshikov, M., Williams, R.J.  Passage-time moments for continuous
non-negative stochastic processes and applications, Adv. Appl. Prob. 1996, 28,
747-762.


\bibitem{Ve1} Veretennikov, A.Yu., On polynomial mixing bounds for stochastic
differential equations. Stochastic Processes and their Applications, 1997, 70,
115-127.

\bibitem{Ve2} Veretennikov, A.Yu. On polynomial mixing and convergence rate for
stochastic difference and differential equations, Teoriya Veroyatnostej i ee
Primenen., 1999, 44(2), 312-327.

\bibitem{Ver01} 
Veretennikov, A.Yu. 
On Polynomial Mixing for SDEs with a Gradient-Type Drift, 
Theory Probab. Appl., 2001, 45(1), 160–164.

\bibitem{Ver81}
Veretennikov, A.Yu. 
On strong and weak solutions of one - dimensional stochastic equations with boundary conditions. Theory Probab. Appl., 1981,  26(4), 670-686.


\end{thebibliography}
\end{document}